\documentclass[12pt]{article}
\usepackage{amsmath,amsthm,amssymb}

\newtheorem{theorem}{Theorem}[section] 

\newtheorem{cor}{Corollary}[section]
\newtheorem{remark}{Remark}[section]

\newtheorem{proposition}{Proposition}[section]

\numberwithin{equation}{section}

\newcommand{\D}{{\rm d}}
\newcommand{\dx}{\, \D x}

\newcommand{\ds}{\, \D s}
\newcommand{\dr}{\, \D r}

\newcommand{\dis}{\displaystyle}

\newcommand{\msp}{\;\;}
\newcommand{\fsp}{\quad\;}
\newcommand{\psp}{\,}

\newcommand{\rz}{\mathbb{R}}
\newcommand{\nz}{\mathbb{N}}

\newcommand{\eps}{\varepsilon}
\newcommand{\iom}{\int_{\Omega}}

\newcommand{\klauf}{\left(\begin{array}}
\newcommand{\klzu}{\end{array}\right)}

\title{On a class of variational problems with linear growth and radial symmetry}
\author{Michael Bildhauer \& Martin Fuchs}
\date{}

\newcommand{\reff}[1]{(\ref{#1})}

\usepackage{hyperref}

\hypersetup{
   pdfnewwindow=true,
   colorlinks=false,
   linkcolor=black,
   citecolor=green,
   filecolor=magenta,
   urlcolor=cyan,
}

\begin{document}

\parindent0em
\maketitle

\newcommand{\op}[1]{\operatorname{#1}}
\newcommand{\bv}{\op{BV}}
\newcommand{\mub}{\overline{\mu}}

\newcommand{\hypref}[2]{\hyperref[#2]{#1 \ref*{#2}}}
\newcommand{\hypreff}[1]{\hyperref[#1]{(\ref*{#1})}}

\newcommand{\ob}[1]{^{(#1)}}
\newcommand{\omin}{\ob{-}}
\newcommand{\omax}{\ob{+}}

\newcommand{\uh}{\hat{u}}
\newcommand{\ut}{{u^{\Theta}}}
\newcommand{\uht}{{\hat{u}^{\Theta}}}
\newcommand{\utn}{{u^{\Theta_0}}}
\newcommand{\uhtn}{{\hat{u}^{\Theta_0}}}

\begin{abstract}
We discuss variational problems on two-dimensional domains with energy densities of
linear growth and with radially symmetric data.\\
The smoothness of generalized minimizers is established under rather weak ellipticity
assumptions. Further results concern the radial symmetry of solutions as well as a precise
description of their behavior near the boundary.\footnote{AMS-Classification: 49J45, 49N60}
\end{abstract}

\section{Introduction}\label{intro}
Inspired by the fundamental work of Giaquinta, Modica and Sou\v{c}ek (\cite{GMS:1979_1}, \cite{GMS:1979_2})
we here discuss the particular minimization problem
\begin{equation}\label{intro 1}
J[w] := \iom g\big(|\nabla w|\big)\dx \to \min \fsp\mbox{in}\msp u_0 + W^{1,1}_0(\Omega)\psp ,
\end{equation}
where $\Omega \subset \rz^2$ is the annulus $\big\{x\in \rz^2:\psp \rho_1 < |x| <\rho_2\big\}$
with radii $0 < \rho_1 < \rho_2 < \infty$.
The function $u_0$ is radially symmetric, which means
\begin{equation}\label{intro 2}
u_0(x) = \hat{u}_0(|x|)\psp , \fsp m_i := \hat{u}_0(\rho_i)\psp ,\msp i=1,\psp 2\psp ,
\end{equation}
reflecting the fact that we want to minimize the functional $J$ among functions with constant values 
on the circles $|x_i| = \rho_i$, $i=1$, $2$. Moreover, we assume that $g\big(|\nabla u|\big)$
is of linear growth w.r.t.~$|\nabla u|$.\\

The purpose of our note is threefold.
\begin{enumerate}
\item We give a general regularity theory for radially symmetric solutions.
In particular, we exclude the occurrence of
an autonomous counterpart of the famous singular example of 
Giaquinta, Modica and Sou\v{c}ek (see \cite{GMS:1979_2}, see also the twodimensional variant given in 
\cite{Bi:1818}). 
Note that we establish the smoothness of solutions up to the boundary which essentially differs from the
attainment of the boundary data (compare \reff{intro 15}). 

\item We allow a wide range of ellipticity since we do not require a balancing condition like 
\begin{equation}\label{balancing}
\frac{g''(s)}{g''(t)}\leq C \fsp\mbox{for all}\msp  s \geq 1\msp\mbox{and}\msp t\in[s/2,2s] \psp ,
\end{equation}
which is a part of the main assumption in \cite{BBM:2018_1}. In fact, this condition is used
for the construction of barriers such that the attainment 
of boundary data can be proved as in \cite{BBM:2018_1} supposing (1.9) of that paper.

Since our arguments leading to the regularity of solutions do not incorporate some detailed estimates
concerning the first derivative of the energy density, we also do not impose an analogue to
(1.7) of Theorem 1.1 given in \cite{BMSV:2018_1}.

\item Following \cite{BF:2018_1}, it is easily shown that boundary data ar respected at least
for $|x|=\rho_2$ which gives the uniqueness of solutions.

Moreover, the possible non-attainment of the boundary data 
in the radially symmetric case has a complete interpretation.\\

\end{enumerate}

Of course we first have to introduce the problem more precisely.\\

In what follows $g$: $[0,\infty) \to [0,\infty)$ is a function of class $C^2\big([0,\infty)\big)$ satisfying
(with suitable constants $a$, $A > 0$, $b$, $B \in \rz$)
\begin{equation}\label{intro 3}
at - b \leq g(t) \leq A t + B \fsp\mbox{for all}\msp t \geq 0
\end{equation}
as well as
\begin{equation}\label{intro 4}
0 = g(0) = g'(0) \psp .
\end{equation}

Let us require for the moment that we just have the inequality
\begin{equation}\label{intro 5}
\nu_1 (1+t)^{-\mu} \leq g''(t) \leq \nu_2 (1+t)^{-1}\psp , \fsp t\geq 0\psp , 
\end{equation}
for some exponent $\mu > 1$, $\nu_1$, $\nu_2$ denoting positive constants. This type of 
$\mu$-elliptic integrand occurs as a special case of the densities discussed, for instance, in \cite{BF:2003_1}
and a series of further papers.\\

Observe that the minimal surface case is included with the choice $g(t) = \sqrt{1+t^2}-1$
leading to \reff{intro 5} with $\mu =3$. Other examples are ($\mu >1$, $k>1$, compare Section \ref{ex})
\begin{eqnarray}
\Phi_\mu(t) &:=& (\mu-1) \int_0^t\int_0^s (1+r)^{-\mu}\dr\ds\psp , \label{intro 6}\nonumber\\[2ex]
 &=& \left\{\begin{array}{rcl}
\dis t - \frac{1}{2-\mu} (1 + t)^{2- \mu} -  \frac{1}{\mu - 2}&\mbox{if}& \mu\not=2\psp ,\\[2ex]
\dis t - \ln(1+t)&\mbox{if}& \mu =2 \psp ,
\end{array}\right.\\[2ex]
\tilde{g}_k(t) &:=& (1+t^k)^{\frac{1}{k}} - 1\psp ,\fsp t \geq 0 \psp ,\label{intro 7}
\end{eqnarray}
where in the latter case we have \reff{intro 5} with $\mu = k+1$. We also note that in recent years 
the example from \reff{intro 7} becomes more and more popular and
in some sense serves as a model for strain-limiting elastic models with linear growth
(see, for instance, \cite{BBMS:2017_1}, \cite{BMS:2015_1},
\cite{BMRS:2014_1} and \cite{BMRW:2015_1}).\\

Associated to our density is the strictly convex integrand
\[
G: \psp \rz^2\to[0,\infty)\psp ,\fsp G(p) := g(|p|) \psp ,\fsp p \in \rz^2 \psp ,
\]
being of linear growth and satisfying the common condition of $\mu$-ellipticity
\begin{equation}\label{intro 8}
\nu_1 \big(1+|p|\big)^{-\mu}|q|^2 \leq D^2G(p)(q,q) \leq \nu_2 (1+|p|)^{-1}|q|^2\psp .
\end{equation}

In fact, \reff{intro 8} follows from the formula
\begin{equation}\label{intro 9}
D^2G(p)(q,q)= \frac{1}{p} g'\big(|p|\big) \Bigg[|q|^2 - \frac{(p\cdot q)^2}{|p|^2}\Bigg]
+g''\big(|p|\big) \frac{(p\cdot q)^2}{|p|^2}
\end{equation}
in combination with \reff{intro 5}.\\

Let us return to our variational problem \reff{intro 1}. As a matter of fact, the existence
of a solution in the subclass $u_0 + W^{1,1}_0(\Omega)$ of the non-reflexive Sobolev space
$W^{1,1}(\Omega)$ (see, e.g., \cite{Ad:1975_1} for a definition of the various Sobolev classes
$W^{k,p}(\Omega)$ and their local variants) can not be guaranteed.
Therefore one has to pass to a suitable relaxed version of \reff{intro 1}.
This approach to linear growth problems is nowadays standard and outlined, for example,
in the monographs \cite{AFP:2000_1}, \cite{Gi:1984_1} and  \cite{GMS:1998_1}, \cite{GMS:1998_2}.
A comprehensive survey of the topic including the historical background is also 
presented in the more recent paper \cite{BS:2013_1}.\\

A relaxed version of \reff{intro 1} is given by
\begin{eqnarray}\label{intro 10}
K[w] &:=& \iom G(\nabla^a w) \dx + \iom G_\infty\Bigg(\frac{\nabla^s w}{|\nabla^s w|}\Bigg)
\,\D|\nabla^s w|\nonumber\\
&&+ \int_{\partial \Omega} G_\infty\Big((u_0-w)\mathcal{N}\Big) \,\D\mathcal{H}^1
\to \min\fsp\mbox{in}\msp \bv(\Omega)\psp , 
\end{eqnarray} 
where $\mathcal{N}$ is the outward unit normal to $\partial \Omega$, $G_\infty$
is the recession function of $G$, and $\nabla^a w$, $\nabla^s w$ denote the regular and
the singular part of $\nabla w$ w.r.t.~the Lebesgue measure. For a definition of the space
$\bv(\Omega)$ we refer to \cite{AFP:2000_1} or \cite{Gi:1984_1}.\\ 

From the convexity of $G$ together with the linear growth condition we obtain
the boundedness of $DG$, moreover, 
\[
g'_\infty := \lim_{t\to \infty} g'(t) = \lim_{t\to \infty} \frac{g(t)}{t}
\]
exists in $(0,\infty)$ and the recession function is given by
\[
G_\infty(p) = g'_\infty |p|\psp ,\fsp p \in \rz^2\psp .
\]
Thus \reff{intro 10} simply reads
\begin{eqnarray}\label{intro 11}
K[w] &=& \iom g\big(|\nabla^a w|\big) \dx + g'_\infty |\nabla^s w|(\Omega)\nonumber\\
&&+ g'_\infty \int_{\partial \Omega} |u_0-w| \,\D\mathcal{H}^1
\to \min\fsp\mbox{in}\msp \bv(\Omega)\psp , 
\end{eqnarray} 
and clearly it holds
\[
K[w] = J[w]\psp ,\fsp \mbox{whenever}\fsp w\in u_0+W^{1,1}_0(\Omega)\psp .
\]

We summarize some known results in the following proposition.

\begin{proposition}\label{intro prop 1}
Let the conditions \reff{intro 2}, \reff{intro 3}-\reff{intro 5} hold for some exponent \mbox{$\mu > 1$}.
Then we have:
\begin{enumerate}
\item The functional $K$ is lower semicontinuous w.r.t.~convergence in $L^1(\Omega)$.
\item Problem \reff{intro 11} admits at least one solution $u \in \bv(\Omega)$.
\item It holds:
\[
\inf_{u_0 +W^{1,1}_0(\Omega)} J = \inf_{\bv(\Omega)} K \psp .
\]
\item $u\in \bv(\Omega)$ is $K$-minimizing $\Leftrightarrow$ $u \in \mathcal{M}$,
\begin{eqnarray*}
\mathcal{M} &:=& \Big\{ v\in L^1(\Omega):\psp \mbox{v is a $L^1(\Omega)$-cluster point}\\
&& \quad\mbox{of some $J$-minimizing sequence from $u_0+W^{1,1}_0(\Omega)$}\Big\} \psp .
\end{eqnarray*}
\item Suppose that \reff{intro 11} admits a solution 
$u\in \bv(\Omega) \cap C^{1}(\Omega)$.
Then any solution $v$ is of the form $v = u+c$ for some $c\in \rz$. Moreover, it
holds $u(x) = \hat{u}(|x|)$.
\item For any $K$-minimizer $u$ we have
\[
\sup_\Omega |u| \leq \max\big\{|m_1|, |m_2|\big\} \psp .
\]
\end{enumerate}
\end{proposition}

In fact, the proposition is based on classical results as the representation formula
of Goffman and Serrin (\cite{GS:1964_1}) and Rehetnyak's continuity theorem
(\cite{Re:1968_1}). We refer to \cite{Bi:1818}, Appendix A, for a detailed discussion
of $iii)$ and $iv)$ which in particular leads to uniqueness Theorem A.9 stated there,
hence to $v)$. Note that a variant of the mentioned Theorem A.9 is also given in \cite{BS:2015_1}, 
Corollary 2.5.
Finally, the last claim is due to Corollary 1 of \cite{BF:2011_1}.\\

Part $v)$ of Proposition \ref{intro prop 1} raises the first challenging question under
which conditions a regular solution 
$u \in \bv(\Omega) \cap C^{1}(\Omega) \subset W^{1,1}(\Omega)$
exists which is immediately leading to the second question, if this minimizer
takes the boundary values $u_0$ thereby solving \reff{intro 1}.\\

Roughly speaking, we have a positive answer to the first problem provided that
\begin{equation}\label{intro 12}
\mu < 3
\end{equation}
(see, e.g., \cite{Bi:1818} or \cite{BS:2013_1}), and from the work of Beck, Bul\'{i}\v{c}ek
and Maringov\'{a} \cite{BBM:2018_1} we deduce that $u=u_0$ on $\partial \Omega$, if \reff{intro 12}
is replaced by the requirement $\mu < 2$ and if the second inequality in \reff{intro 5}
is replaced by the condition $g''(t) \leq \nu_2 (1+t)^{-\mu}$.\\

In the situation at hand we neither require
any upper bound on the exponent $\mu$ nor a balancing condition in the sense of \reff{balancing}
still giving a positive answer to the existence of a smooth $K$-minimizer.\\

We just need the limitation \reff{intro 14}
for the range of anisotropy admissible in the behavior of $g''$, which is quite similar
to the superlinear analogue $q< p+2$ in the case of anisotropic growth conditions
(see \cite{Bi:1818} for an overview and a list of references).\\

Let us now state our main results.

\begin{theorem}\label{intro theo 1}
Suppose that $\mu \in (1,\infty)$, let \reff{intro 2}, \reff{intro 3}, \reff{intro 4} hold and replace
\reff{intro 5} by the condition
\begin{equation}\label{intro 13}
\nu_1 (1+t)^{-\mu} \leq g''(t) \leq v_2 (1+t)^{-\mub}\psp , \fsp t \geq 0 \psp ,
\end{equation}
for some exponent $\mub \in [1,\mu]$ such that
\begin{equation}\label{intro 14}
\mu -\mub < 2 \psp .
\end{equation}
Then the relaxed problem \reff{intro 11} admits a solution
\[
u \in W^{1,1}(\Omega)\cap C^0(\overline{\Omega}) \cap C^1(\Omega)\cap W^{2,2}_{\op{loc}}(\Omega)
\]
which in addition is of the form $u(x) = \hat{u}\big(|x|\big)$. Moreover, the solution is unique up to
additive constants.
\end{theorem}

\begin{remark}\label{intro ex 1}
\begin{enumerate} 
\item In the case $\mu$, $\mub > 1$ we deduce from \reff{intro 13} and \reff{intro 4}
the inequality 
\[
c \Phi_{\mu}(t) \leq g(t) \leq C \Phi_{\mub} (t)
\]
with $\Phi_{\dots}$ defined in \reff{intro 6}, which means that \reff{intro 3}
automatically holds.

\item Note that in particular the one parameter family of energy densities given in \reff{intro 7} and suitable
generalizations are covered by our considerations.

\item We may take any function $\psi(t)$ satisfying for some constants
$c_1$, $c_2 \in \rz$
\[
c_1 (1+ t)^{-\mu} \leq \psi(t) \leq c_2 (1+t)^{-\mub}\psp ,\fsp t\geq 0 \psp ,
\]
and obtain a function
\[
\Psi(t) = \int_0^t \int_0^s \psi(r) \dr \ds 
\]
which clearly satisfies \reff{intro 13} but in general violates a balancing condition like given
in (1.8) of \cite{BBM:2018_1}.

\end{enumerate}
\end{remark}

We do not know if the solution $u$ takes the boundary values $u_0$
for $|x| = \rho_1$.
However, the following theorem yields a complete description of the boundary behavior.

\begin{theorem}\label{intro theo 2}
The minimizer given in Theorem \ref{intro theo 1} in fact is the unique solution of problem
\reff{intro 11}. Moreover, we have:
\begin{enumerate}
\item The minimizer respects the boundary data for $|x| = \rho_2$, hence it is the solution
of the minimizing problem
\begin{eqnarray}\label{intro 15}
\iom g\big(|\nabla w|\big)\dx  + g'_\infty \int_{|x| = \rho_1} |w-m_1| \, \D\mathcal{H}^1 \to 
\min&&\nonumber \\[1ex]
w \in W^{1,1}(\Omega) \psp ,\fsp w=m_2 \msp\mbox{on} \msp\big\{|x| = \rho_2\big\}\psp .&&
\end{eqnarray}
\item Suppose that $m_2$ is fixed, abbreviate $m=m_1$ and 
let $u_m(x) = \uh_m\big(|x|\big)$ denote the unique solution of \reff{intro 15}.

Suppose w.l.o.g.~that $m < m_2$. Then we have
\begin{enumerate}
\item For any $\rho\in (\rho_1, \rho_2)$ it holds that $\hat{u}_m(\rho) \geq m$.
\item As a function depending on $m$, the quantity $\hat{u}_m(\rho_1)$ is a non-decreasing function, i.e.
\[
\zeta_1 < \zeta_2\fsp\Rightarrow \fsp \uh_{\zeta_1}(\rho_1) \leq \uh_{\zeta_2}(\rho_2) \psp .
\]
\end{enumerate}
\end{enumerate}
\end{theorem}

As a corollary we obtain in particular: 

\begin{cor}\label{intro cor 1}
With the notation of Theorem \ref{intro theo 2} we suppose that there exists $m< m_2$ such that
the boundary data are not attained for $|x|=\rho_1$.

Then for any $\zeta \leq m$ we have $u_{\zeta} \equiv u_{m}$. 
\end{cor}

\section{Proof of Theorem \ref{intro theo 1}}\label{proof}

We proceed by induction showing that the statements of the theorem hold provided
$\mu \in (1,k)$ for some $k\geq 2$.\\

Let in the beginning $k=3$. From \reff{intro 13} we immediately get \reff{intro 8}
(recall \reff{intro 9}) and from Theorem 4.32 in \cite{Bi:1818} we deduce the existence
of a unique (up to constants) generalized minimzer $u$ of class 
$\bv(\Omega) \cap C^1(\Omega) \subset W^{1,1}(\Omega)$.
Alternatively, we can quote Theorem 4.16 from this reference observing that
Assumption 4.11 trivially holds for the situation at hand.\\

Proposition \ref{intro prop 1}, $v)$, implies that $u(x) = \hat{u}\big(|x|\big)$ with
$\hat{u} \in W^{1,1}(\rho_1,\rho_2) \subset C^0\big([\rho_1,\rho_2]\big)$ 
(see \cite{BGH:1}, Chapter 2), hence $u \in C^0(\overline{\Omega})$. In order to
show
\begin{equation}\label{proof 1}
u \in W^{2,2}_{\op{loc}}(\Omega)
\end{equation}
it is sufficient to prove uniform local $W^{2,2}$-bounds for the solutions $u_\delta$ of
the regularized problem
\[
J_\delta[w]:= \frac{\delta}{2} \iom |\nabla w|^2 \dx + \iom g\big(|\nabla w|\big)\dx
\to \min \fsp\mbox{in}\msp u_0+ W^{1,2}_0(\Omega) \psp .
\]
To this purpose we just quote Lemma 4.19, $i)$, in \cite{Bi:1818} choosing the exponent $s$
so large that the l.h.s.~of the Caccioppoli inequality is bounded from below by
\[
\alpha \iom \eta^2 |\nabla^2 u_\delta|^2 \dx \psp ,\fsp\mbox{$\alpha$ denoting a suitable
uniform constant}\psp .
\]
On the r.h.s.~we observe Theorem 4.25 from \cite{Bi:1818}, which gives the desired
uniform bound for $u_\delta \in W^{2,2}_{\op{loc}}(\Omega)$ leading to \reff{proof 1}.\\

Suppose next that $k\geq 3$ and that Theorem \ref{intro theo 1} is true for exponents
$\mu \in (1,k)$. We then claim the validity of Theorem \ref{intro theo 1} for
\begin{equation}\label{proof 2}
\mu \in [k,k+1) \psp .
\end{equation}
So let the density $g$ satisfy \reff{intro 13} with exponent $\mub \in [1, \mu]$
such that \reff{intro 14} is true. We choose
\begin{equation}\label{proof 3}
\tau \in \big(\mu -1, \min\{k,\mub\}\big)
\end{equation}
and observe the inequalities
\begin{equation}\label{proof 4}
\tau < k \psp ,\fsp \tau < \mub\psp , \fsp \mu - \tau < 1\psp .
\end{equation}
For $\delta \in (0,1)$ we introduce the density
\begin{equation}\label{proof 5}
g_\delta(t) := \delta \Phi_\tau(t) + g(t)\psp ,\fsp t \geq 0 \psp ,
\end{equation}
with function $\Phi_\tau$ from \reff{intro 6}. Moreover, we let
\begin{equation}\label{proof 6}
G_\delta (p) := g_\delta\big(|p|\big)\psp , \fsp p \in \rz^2\psp .
\end{equation}
Then it holds 
\[
g''_\delta(t) = \delta (\tau-1)  (1+t)^{-\tau} + g''(t)
\]
and the second inequality in \reff{proof 4} together with \reff{intro 13} shows
\begin{equation}\label{proof 7}
c_1(\delta) (1+t)^{-\tau} \leq g''_{\delta}(t) \leq c_2(\delta) (1+t)^{-\tau}
\end{equation}
with constants $c_i(\delta) > 0$. Recalling the first inequality in \reff{proof 4}
and observing \reff{proof 7} our inductive hypothesis applies to the regularized
problem
\begin{equation}\label{proof 8}
K_\delta [w] \to \min \fsp\mbox{in}\msp \bv(\Omega) \psp ,
\end{equation}
where $K_\delta$ is defined according to \reff{intro 10} and \reff{intro 11} with $G$
and $g$ replaced by $G_\delta$ and $g_\delta$, respectively. Let
\begin{equation}\label{proof 9}
u_\delta \in W^{1,1}(\Omega) \cap C^0(\overline{\Omega}) \cap C^1(\Omega) 
\cap W^{2,2}_{\op{loc}}(\Omega)
\end{equation}
denote the unique (up to constants) solution to \reff{proof 8} which additionally
satisfies $u_\delta(x) = \hat{u}_\delta\big(|x|\big)$. The regularity
properties of $u_\delta$ stated in \reff{proof 9} are in turn sufficient to derive
the Caccioppoli inequality from Lemma 4.19, $ii)$, in \cite{Bi:1818}, i.e. it holds
\begin{eqnarray}\label{proof 10}
\lefteqn{\iom D^2G_\delta(\nabla u_\delta)(\partial_\gamma \nabla u_\delta,
\partial_\gamma \nabla u_\delta) \Gamma_\delta^s \eta^{2l}\dx}\nonumber\\
&\leq& c \iom D^2G_\delta(\nabla u_\delta)(\nabla \eta,\nabla \eta)\eta^{2l-2}
|\nabla u_\delta|^2 \Gamma_\delta^s\dx
\end{eqnarray}
for any $s\geq 0$, $l\in \nz$ and $\eta \in C^1_0(\Omega)$, $0\leq \eta \leq 1$,
where we have abbreviated $\Gamma_\delta := 1+|\nabla u_\delta|^2$ and the sum is taken
w.r.t.~the index $\gamma$.\\

Letting $\eta(x) = \hat{\eta}\big(|x|\big)$ we set
\[
p= \uh'_\delta\big(|x|\big) \frac{x}{|x|} \psp ,\fsp q = \hat{\eta}'\big(|x|\big) \frac{x}{|x|} \psp ,
\]
and observe that we have in \reff{intro 9}
\[
\Bigg[ |q|^2 - \frac{(p\cdot q)^2}{|p|^2}\Bigg] = 0 \psp .
\]

This reduces \reff{proof 10} to the inequality
\begin{equation}\label{proof 11}
\iom g''_\delta \big(|\nabla u_\delta|\big)|\nabla^2 u_\delta|^2 \eta^{2l} \Gamma_\delta^s \dx
\leq c \iom g''_\delta\big(|\nabla u_\delta|\big) \Gamma_\delta^{s+1} |\nabla \eta|^2 \eta^{2l-2}\dx
\end{equation}
with constant $c>0$ not depending on $\delta$. Applying \reff{intro 13} and recalling the definition 
\reff{proof 5} we arrive at (neglecting the $\delta$-term on the l.h.s~of \reff{proof 11})
\begin{eqnarray}\label{proof 12}
\lefteqn{\iom \eta^{2l} |\nabla^2 u_\delta|^2 \Gamma_\delta^{s-\frac{\mu}{2}} \dx}\nonumber\\
&\leq & c\Bigg[ \delta \iom \eta^{2l-2} |\nabla \eta|^2 \Gamma_\delta^{s+1-\frac{\tau}{2}}\dx
+ \iom \eta^{2l-2} |\nabla \eta|^2 \Gamma_\delta^{s+1-\frac{\mub}{2}}\dx\Bigg] \psp .
\end{eqnarray}

Next we choose $\varphi := u_\delta \Gamma_\delta^{\frac{\alpha}{2}} \eta^{2l}$
as admissible (recall \reff{proof 9}) test function in the Euler equation
\begin{equation}\label{proof 13}
0 = \iom DG_\delta(\nabla u_\delta) \cdot \nabla \varphi \dx \psp ,
\end{equation}
where $\eta$ and $l$ are as above and $\alpha \geq 0$ is some number to be fixed later.
With this choice \reff{proof 13} gives
\begin{eqnarray}\label{proof 14}
\lefteqn{\iom DG_\delta (\nabla u_\delta) \cdot \nabla u_\delta \Gamma_\delta^{\frac{\alpha}{2}} 
\eta^{2l}\dx}\nonumber\\
&=& - \iom DG_\delta(\nabla u_\delta)\cdot \nabla \eta \eta^{2l-1} 2 l u_\delta
\Gamma_\delta^{\frac{\alpha}{2}} \dx\nonumber\\
&& - \iom DG_\delta (\nabla u_\delta)\cdot \nabla \Gamma_\delta^{\frac{\alpha}{2}} \eta^{2l}
u_\delta \dx =: T_1 + T_2 \psp .
\end{eqnarray}
We have
\[
DG_\delta (\nabla u_\delta) \cdot \nabla u_\delta = g'_\delta (t) t \geq g'(t)t \psp ,
\fsp t:= |\nabla u_\delta|\psp ,
\]
and from the first inequality in \reff{intro 13} we get
\[
g'(t)  \geq c \int_0^t (1+s)^{-\mu}\ds \geq c \Big[1-(1+t)^{1-\mu}\Big]
\]
and in conclusion
\[
DG_\delta(\nabla u_\delta) \cdot \nabla u_\delta \geq c \big[|\nabla u_\delta| - 1 \big]\psp ,
\]
where as usual the value of $c$ may vary from line to line. Therefore we get
\begin{equation}\label{proof 15}
\mbox{l.h.s.~of \reff{proof 14}}\geq
c \Bigg[ \iom \Gamma_\delta^{\frac{\alpha +1}{2}}\eta^{2l} \dx 
- \iom \eta^{2l} \Gamma_\delta^{\frac{\alpha}{2}} \dx \Bigg] \psp .
\end{equation}
For $T_1$, $T_2$ on the r.h.s.~of \reff{proof 14} we use (see Proposition \ref{intro prop 1}, $vi$))
\[
\sup_\Omega |u_\delta| \leq \max\big\{|m_1|,|m_2|\big\} \psp ,
\]
as well as the uniform boundedness of 
\[
DG_\delta (p) = g_\delta'(|p|) \frac{p}{|p|} \psp ,
\]
which is immediate by the definition of $g_\delta$ and the properties of $g$. We obtain
\begin{eqnarray*}
|T_1| &\leq & c \iom |\nabla \eta| \eta^{2l-1} \Gamma_\delta^{\frac{\alpha}{2}}\dx \psp ,\\
|T_2| &\leq & c \iom \eta^{2l} |\nabla^2 u_\delta| \Gamma_\delta^{\frac{\alpha -1}{2}} \dx \psp .
\end{eqnarray*}
Returning to \reff{proof 14}, using \reff{proof 15} and the inequalities for $T_i$, it is
shown ($c=c(l)$) that
\begin{eqnarray}\label{proof 16}
\iom \eta^{2l} \Gamma_\delta^{\frac{\alpha +1}{2}} \dx &\leq & 
c \Bigg[ \iom \eta^{2l} \Gamma_\delta^{\frac{\alpha}{2}} \dx 
+ \iom |\nabla \eta| \eta^{2l-1} \Gamma_\delta^{\frac{\alpha}{2}} \dx\nonumber\\
&&+ \iom \eta^{2l} |\nabla^2 u_\delta| \Gamma_\delta^{\frac{\alpha -1}{2}}\dx \Bigg]\nonumber\\
&=:& c \big[S_1 + S_2 + S_3\big]\psp .
\end{eqnarray}
To the quantities $S_i$, $i=1$, $2$, $3$, we apply Young' inequality:
\begin{eqnarray*}
S_1 &\leq& \eps \iom \eta^{2l} \Gamma_\delta^{\frac{\alpha+1}{2}} \dx +
c(\eps) \iom \eta^{2l} \Gamma_\delta^{\frac{\alpha -1}{2}} \dx \psp ,\\[1ex]
S_2 &\leq & \eps \iom \eta^{2l} \Gamma_\delta^{\frac{\alpha+1}{2}} \dx
+c(\eps) \iom \eta^{2l-2} |\nabla \eta|^2 \Gamma_\delta^{\frac{\alpha -1}{2}} \dx \psp ,\\[1ex]
S_3 &\leq & \eps \iom \eta^{2l} \Gamma_\delta^{\frac{\alpha +1}{2}} \dx
+ c(\eps) \iom \eta^{2l} |\nabla^2 u_\delta|^2 \Gamma_\delta^{\frac{\alpha -3}{2}} \dx\psp .
\end{eqnarray*}
For $\eps$ sufficiently small we obtain from \reff{proof 16}
\begin{eqnarray}\label{proof 17}
\lefteqn{\iom \eta^{2l} \Gamma_\delta^{\frac{\alpha +1}{2}} \dx}\nonumber\\
&\leq & c \Bigg[ \iom \eta^{2l} \Gamma_\delta^{\frac{\alpha -3}{2}} |\nabla^2 u_\delta|^2 \dx
+ \iom \eta^{2l-2} \big[\eta^2 + |\nabla \eta|^2\big] \Gamma_\delta^{\frac{\alpha -1}{2}} \dx
\Bigg]\psp .
\end{eqnarray}
In a final step we estimate
\begin{eqnarray*}
\iom \eta^{2l-2} |\nabla \eta|^2 \Gamma_\delta^{\frac{\alpha -1}{2}}\dx
&=& \iom \eta^{2l-2} \Gamma_\delta^{\frac{\alpha +1}{4}} |\nabla \eta|^2 \Gamma_\delta^{\frac{\alpha -3}{4}}\dx\\[1ex]
&\leq & \eps \iom \eta^{4l-4}\Gamma_\delta^{\frac{\alpha +1}{2}} \dx + c(\eps)
\iom |\nabla \eta|^4 \Gamma_\delta^{\frac{\alpha -3}{2}} \dx\\[1ex]
&\leq & \eps \iom \eta^{2l} \Gamma_\delta^{\frac{\alpha +1}{2}} \dx + c(\eps)
\iom |\nabla \eta|^4 \Gamma_\delta^{\frac{\alpha-3}{2}}\dx \psp ,
\end{eqnarray*}
where we have used $\eta^{4l-4} \leq \eta^{2l}$, $l\geq 2$, on account of $0 \leq \eta \leq 1$.
Clearly it holds
\[
\iom \eta^{2l} \Gamma_\delta^{\frac{\alpha -1}{2}} \dx \leq 
\eps \iom \eta^{2l} \Gamma_\delta^{\frac{\alpha+1}{2}} \dx +
c(\eps) \iom \eta^{2l} \Gamma_\delta^{\frac{\alpha -3}{2}}\dx \psp ,
\]
and \reff{proof 17} implies 
\begin{equation}\label{proof 18}
\iom \eta^{2l} \Gamma_\delta^{\frac{\alpha +1}{2}} \dx \leq c \Bigg[
\iom \eta^{2l} \Gamma_\delta^{\frac{\alpha -3}{2}} |\nabla^2 u_\delta|^2\dx
+ \iom \Big[|\nabla \eta|^4 + \eta^{2l}\Big] \Gamma_\delta^{\frac{\alpha -3}{2}} \dx \Bigg]\psp .
\end{equation}
Let us choose $\alpha =3$ in \reff{proof 18} yielding
\begin{equation}\label{proof 19}
\iom \eta^{2l} \Gamma_\delta^2 \dx \leq c \Bigg[ \iom \eta^{2l} |\nabla^2 u_\delta|^2 \dx + c(\eta) \Bigg]\psp .
\end{equation}
On the r.h.s.~of \reff{proof 19} we apply \reff{proof 12} for the choice $s=\mu/2$ yielding
\begin{eqnarray}\label{proof 20}
\iom \eta^{2l} \Gamma_\delta^2 \dx &\leq & c \Bigg[ \delta \iom \eta^{2l-2} |\nabla \eta|^2
\Gamma_\delta^{\frac{\mu}{2}+1-\frac{\tau}{2}} \dx\nonumber\\
&&+ \iom \eta^{2l-2} |\nabla \eta|^2 \Gamma_\delta ^{\frac{\mu}{2} +1 - \frac{\mub}{2}}\dx
+ c(\eta) \Bigg]\psp .
\end{eqnarray}
\reff{proof 4} implies
\[
\frac{\mu}{2} - \frac{\tau}{2} + 1 < \frac{3}{2}
\]
and from \reff{intro 14} it follows
\[
\frac{\mu}{2} - \frac{\mub}{2}+1 < 2 \psp ,
\]
thus we have to handle terms like
\[
\iom \eta^{2l-2} \Gamma_\delta^p |\nabla \eta|^2 \dx 
\]
with exponent $p\in (1,2)$ on the r.h.s.~of \reff{proof 20}.
Evidently it holds for $l$ sufficiently large
\[
\iom \eta^{2l-2} |\nabla \eta|^2 \Gamma_\delta^p \dx \leq \eps \iom \eta^{2l} \Gamma_\delta^2 \dx 
+ c(\eps,\eta)
\]
and therefore \reff{proof 20} implies
\begin{equation}\label{proof 21}
|\nabla u_\delta| \in L^4_{\op{loc}}(\Omega)\fsp\mbox{uniformly in}\msp \delta \psp .
\end{equation}
Next we let $\alpha =7$ in \reff{proof 18} and $s=2 +\mu/2$ in \reff{proof 12}.
Taking into account \reff{proof 21} a repetition of the preceeding calculations
leads to $|\nabla u_\delta| \in L^8_{\op{loc}}(\Omega)$ and by iteration we find
for any $q < \infty$
\begin{equation}\label{proof 22}
|\nabla u_\delta| \in L^q_{\op{loc}}(\Omega)\fsp\mbox{uniformly in}\msp \delta \psp .
\end{equation}

With \reff{proof 22} we deduce from \reff{proof 12} with the choice $s=\mu/2$ uniform higher
weak differentiability, i.e.
\begin{equation}\label{proof 23}
u_\delta \in W^{2,2}_{\op{loc}}(\Omega) \fsp\mbox{uniformly in}\msp \delta\psp .
\end{equation}
From
\[
K[u_\delta] \leq K_\delta[u_\delta] \leq K_\delta[u_0] \leq c(u_0) < \infty
\]
together with
\[
\|u_\delta\|_{L^\infty(\Omega)} \leq \max\big\{|m_1|,|m_2|\big\}
\]
it follows
\[
\sup_{0 < \delta < 1} \|u_\delta\|_{W^{1,1}(\Omega)} < \infty \psp ,
\]
hence there is a function $u \in \bv(\Omega)$ such that
\begin{equation}\label{proof 24}
u_\delta \to u\fsp\mbox{in}\msp L^1(\Omega)
\end{equation}
at least for a subsequence. We claim that $u$ is $K$-minimizing.
Let $v \in \bv(\Omega)$. By Proposition \ref{intro prop 1}, $i)$, and \reff{proof 24}
it holds
\[
K[u] \leq \lim_{\delta \to 0} K[u_\delta] \psp .
\]
At the same time we have by the minimizing property of $u_\delta$
\[
K[u_\delta] \leq K_\delta[u_\delta] \leq K_\delta[v]
\to K[v] \fsp\mbox{as}\msp \delta \to 0 \psp ,
\]
which proves our claim. Obviously \reff{proof 24} implies the validity of
\reff{proof 22} and \reff{proof 23} for the function $u$. Moreover,
the radial symmetry of $u_\delta$ extends to $u$. This completes the
proof of Theorem \ref{intro theo 1}. \hfill \qed

\section{Proof of Theorem \ref{intro theo 2} and Corollary \ref{intro cor 1}}\label{geo}
\underline{\it Proof of Theorem \ref{intro theo 2}}\\

Suppose that $\tilde{u}$ is any given solution of \reff{intro 11}. 
The first part of Theorem \ref{intro theo 1} guarantees that $\tilde{u}$ is sufficiently
smooth such that any solution of \reff{intro 11} is of the form $\tilde{u}+c$, $c\in \rz$.\\

In order to show uniqueness together with claim $i)$, we distinguish four different cases:\\

{\it Case 1. The data are attained on the whole boundary $\partial \Omega$.}\\

Then, if $\tilde{u}+c$, $c\in \rz$, is a candidate for a possibly different minimizer, then on  account of 
\[
0 = \int_{\partial \Omega} |\tilde{u}-u_0| \, \D\mathcal{H}^1 
= \int_{\partial \Omega} |\tilde{u}+c-u_0| \, \D\mathcal{H}^1 = |c| \mathcal{H}^1(\partial \Omega)\psp ,
\]
$c=0$ is immediate, hence $\tilde{u}$ is the unique solution. (Case 1 corresponds to \cite{BS:2013_1}, 
Lemma 5.5.)\\

{\it Case 2. Both for $|x|=\rho_1$ and for $|x|= \rho_2$ the solution $\tilde{u}$ does not 
attain the boundary data.}\\

Following \cite{BF:2018_1}, we  let for any $w\in \bv(\Omega)$
\begin{eqnarray*}
\partial_+^w\Omega &:=& \big\{x\in \partial \Omega:\psp w(x) > u_0(x)\big\}\psp ,\\[1ex]
\partial_-^w\Omega &:=& \big\{x\in \partial \Omega:\psp w(x) < u_0(x)\big\}\psp ,\\[1ex]
\partial_0^w\Omega &:=& \big\{x\in \partial \Omega:\psp w(x) = u_0(x)\big\}\psp ,
\end{eqnarray*}

and observe that Theorem 2.4 of this reference just needs the hypothesis
of the strict convexity of the linear growth energy density. Thus, Theorem 2.4
shows for the solution $\tilde{u}$
\begin{equation}\label{uni 1}
\big|\mathcal{H}^1(\partial_+^{\tilde{u}} \Omega) - \mathcal{H}^1(\partial_-^{\tilde{u}} \Omega)\big| \leq 
\mathcal{H}^1(\partial_0^{\tilde{u}} \Omega) \psp .
\end{equation}

Since the boundary data are completely ignored in the case under consideration, we have
\begin{equation}\label{uni 2}
\mathcal{H}^1(\partial_0^{\tilde{u}} \Omega)=0\quad\mbox{and in conclusion}\quad
\mathcal{H}^1(\partial_+^{\tilde{u}} \Omega) = \mathcal{H}^1(\partial_-^{\tilde{u}} \Omega) \psp .
\end{equation}
This, however, is not possible on account of
\begin{equation}\label{uni 3}
\Big|\big\{|x| = \rho_1\big\}\Big| < \Big|\big\{|x| = \rho_2\big\}\Big|\psp .
\end{equation}

{\it Case 3. The boundary data are attained for $|x|= \rho_1$, they are not attained for $|x|=\rho_2$.}\\

In this case
\[
\partial_0^{\tilde{u}} \Omega =\big\{|x|=\rho_1\big\}
\]
gives a contradiction referring to \reff{uni 1} and \reff{uni 3}.\\

{\it Case 4. The boundary data are attained for $|x|= \rho_2$, they are not attained for $|x|=\rho_1$.}\\

This case is possible and in accordance with our claim
\[
\tilde{u} = m_2 \fsp\mbox{on}\msp \big\{|x| = \rho_2\big\}\fsp
\mbox{for any solution $\tilde{u}$ of \reff{intro 11}}\psp .
\]

Since by Theorem \ref{intro theo 1} uniqueness up to additive constants holds true, we now even have 
the uniqueness of solutions on account of the attainment of the data for $|x| = \rho_2$.\\

Next we prove our claim $ii)$. In the following $m_2$ is fixed and we suppose by the first part of the theorem
that for any solution under consideration we have
\[
u(\rho_2) = u_0(\rho_2) = m_2 \fsp\mbox{and w.l.o.g.}\fsp m_2> m_1 =: m\psp .
\]
In the case $m_2 < m_1$ the analogous arguments are obvious.\\

Let us define for any $w\in \bv(\Omega)$ satisfying $w=m_2$ for $|x|=\rho_2$ and for 
any real number $\zeta < m_2$ the energies
\begin{eqnarray*}
K^0[w] &:=& \iom g \big(|\nabla^a w|\big) \dx \psp ,\\[2ex]
K_\zeta[w] &:=&\iom g\big(|\nabla^a w|\big)\dx  + g'_\infty |\nabla^s w|(\Omega)
+ g'_\infty \int_{|x| = \rho_1} |w-\zeta| \, \D\mathcal{H}^1 \psp .
\end{eqnarray*}
By the first part and by Theorem \ref{intro theo 1}, the unique solution 
$u_\zeta(x) = \hat{u}_\zeta\big(|x|\big)$
of the minimizing problem
\[
K_\zeta[w]\to \min  \fsp \mbox{in}\msp \bv(\Omega)
\] 
in particular is of class $W^{1,1}(\Omega)$, hence
\begin{equation}\label{mono 0}
K_\zeta[u_\zeta] = \iom g\big(|\nabla u_\zeta|\big) \dx + 
g'_\infty \int_{|x| = \rho_1} |u_\zeta- \zeta| \, \D\mathcal{H}^1 
\end{equation}
and $K_\zeta[w]$ takes the form \reff{mono 0} whenever $w \in W^{1,1}(\Omega)$.\\

Establishing our claim (a) we suppose by contradiction that there exists
$\rho \in (\rho_1,\rho_2)$ such that $\hat{u}_\zeta(\rho) < \zeta$.\\ 

Then the continuity of $u_\zeta$ yields a real number 
$\hat{\rho}\in (\rho,\rho_2)$ such that $\uh_{\zeta}(\hat{\rho}) = \zeta$
and the choice
\[
w_\zeta(x) :=\left\{\begin{array}{rcl}
u_\zeta(x) &\mbox{for}& |x|\in (\hat{\rho},\rho_{2})\psp ,\\[1ex]
\zeta &\mbox{for}& |x|\in (\rho_1,\hat{\rho}]\psp 
\end{array}\right.
\]
immediately contradicts the minimality of $u_\zeta$.\\

In order to prove claim (b) we suppose that there exist real numbers 
\begin{equation}\label{mono 1}
\zeta_1 < \zeta_2 \fsp\mbox{and}\fsp 
\uh_{\zeta_2}(\rho_1) = \zeta_2\omax < \zeta_1\omax = \uh_{\zeta_1}(\rho_1) \psp . 
\end{equation}

Part (a) shows that in this case we have 
\begin{equation}\label{mono 2}
\zeta_1 < \zeta_2 \leq \zeta_2\omax < \zeta_1\omax < m_2
\end{equation}
which guarantees the positive sign of the penalty terms below.\\
 
Note that, given two real numbers $\xi$, $\kappa$ such that $m_2 \geq \xi\geq \kappa$, 
part (a) also implies the representation formula
\begin{eqnarray}\label{mono 3}
K_\xi[u_\xi] & =& K^0[u_\xi] + g'_\infty \int_{|x|=\rho_1} |u_\xi -\xi|\, \D\mathcal{H}^{n-1}\nonumber\\[1ex]
&=& K^0[u_\xi] + g'_\infty 2\pi \rho_1 \big(\uh_\xi(\rho_1) -\xi\big)\nonumber \\[2ex] 
&=& K^0[u_\xi] + g'_\infty 2 \pi \rho_1 \big(\hat{u}_\xi(\rho_1)-\kappa\big) - 
g'_\infty 2\pi \rho_1 (\xi-\kappa )\nonumber\\[2ex]
&=& K_\kappa[u_\xi] - g'_\infty 2\pi\rho_1 (\xi-\kappa)  \psp .
\end{eqnarray}

Now we proceed by observing
\begin{eqnarray}\label{mono 4}
K_{\zeta_1} [u_{\zeta_1}] &=& K^0[u_{\zeta_1}] 
+ g'_\infty 2 \pi \rho_1 (\zeta_1\omax -\zeta_1)\nonumber\\[2ex]
&=& K^0[u_{\zeta_1}] + g'_\infty 2 \pi \rho_1 (\zeta_1\omax - \zeta_2)\nonumber \\[2ex]
&& + g'_\infty 2 \pi \rho_1 (\zeta_2-\zeta_1)\nonumber \\[2ex]
&\geq & K_{\zeta_2}[u_{\zeta_2}]+ g'_\infty 2 \pi \rho_1 (\zeta_2-\zeta_1)\nonumber\\[2ex]
&= & K_{\zeta_1}[u_{\zeta_2}]  \psp ,
\end{eqnarray}
where we recall \reff{mono 2} for the discussion of the absolute values in the penalty term. 
Moreover, the inequality
\[
K^0[u_{\zeta_1}] + g_{\infty}' 2 \pi \rho_1 (\zeta_1\omax- \zeta_2) = K_{\zeta_2}[u_{\zeta_1}] 
\geq K_{\zeta_2}[u_{\zeta_2}] 
\]
follows from the minimality of $u_{\zeta_2}$ and the last equality in 
\reff{mono 4} is due to \reff{mono 3}.\\

Finally we observe that inequality \reff{mono 4} would give $u_{\zeta_1} = u_{\zeta_2}$ 
by uniqueness of minimizers which contradicts
the hypothesis \reff{mono 1}, i.e.~we have a contradiction to $\zeta_2\omax < \zeta_1\omax$,
and the proof of Theorem \ref{intro theo 2} 
is complete. \hfill \qed\\

\underline{\it Proof of Corollary \ref{intro cor 1}.}\\

Using the notation of Theorem \ref{intro theo 1} we first recall 
two facts that are already established above:
\begin{enumerate}
\item $\uh \in W^{1,1} (\rho_1,\rho_2) \cap C^1(\rho_1,\rho_2)$;
\item $\uh(\rho_2) = m_2$.
\end{enumerate}

For the reader's convenience we sketch some general observations on the Euler equation
which can also be found in \cite{BBM:2018_1}:\\

Given a test function $\eta\in C^\infty_0 (\Omega)$ we have
\begin{equation}\label{cor 1}
\iom \frac{g'\big(|\nabla u|\big)}{|\nabla u|} \nabla u \cdot \nabla \eta \dx = 0  \psp .
\end{equation}
Inserting
\[
\nabla u = \uh'\big(|x|\big) \frac{x}{|x|} \psp .
\]
and choosing $\eta(x) = \hat{\eta}\big(|x|\big)$ we obtain
\begin{equation}\label{cor 2}
0 = \iom g'\big(|\uh'|\big) \frac{\uh'\big(|x|\big)}{\big|\uh'\big(|x|\big)\big|}
\hat{\eta}'\big(|x|\big) \dx
= 2 \pi \int_{\rho_1}^{\rho_2}  g'\big(|\uh'|\big) \frac{\uh'}{|\uh'|} \hat{\eta}' r \dr \psp . 
\end{equation}

Note that on account of \reff{intro 4} the expression $g'(t)/t$ is well defined
in the limit $t \to 0$.\\

Using \reff{cor 2}, Du Bois-Reymond's lemma as variant of the fundamental lemma
implies the existence of a real number $\lambda \in \rz$ such that
\begin{equation}\label{cor 3}
g'\big(|\uh'|\big) \frac{\uh'}{|\uh'|} = \frac{\lambda}{r} \psp .
\end{equation}

If $\uh \not\equiv 0$, then zeroes of $\uh$ are excluded by \reff{cor 3} and supposing
w.l.o.g.~$m_1 < m_2$ we have $\uh' >0$ and \reff{cor 3} reduces to
\begin{equation}\label{cor 4}
g' \big(\uh'(r)\big) = \frac{\lambda}{r}\fsp\mbox{for all}\msp r\in (\rho_1,\rho_2) \psp .
\end{equation}

By assumption $g$ is a strictly convex function, i.e.~$g'$ is a strictly
increasing function and we have that
\[
g':\psp (0,\infty) \to (0,g_\infty')\fsp\mbox{is one-to-one}\psp ,
\]
hence we obtain from \reff{cor 4}
\begin{equation}\label{cor 5}
0 < \uh'(r) = (g')^{-1} (\lambda/r) \fsp\mbox{for all}\msp r\in (\rho_1,\rho_2) \psp .
\end{equation}
Note the validity of \reff{cor 4} for all $r\in (\rho_1,\rho_2)$
and in conclusion the possible values of $\lambda$ are given by
\begin{equation}\label{cor 6}
0 < \lambda \leq \rho_1 g'_\infty \psp . 
\end{equation}
Finally we consider the possible range for realizing boundary data:
\[
\Delta m (\lambda) := u(\rho_2) - u(\rho_1) = \int_{\rho_1}^{\rho_2} (g')^{-1} (\lambda/r) \dr \psp ,
\]
where we note that $\Delta m (\lambda) \to 0$ as $\lambda \to 0$.\\

Now, on account of \reff{cor 6}, for any $\rho_1 < \hat{\rho} < \rho_2$ the function
$(g')^{-1}(\lambda/r)$ is bounded in $(\hat{\rho},\rho_2]$ with a constant not depending
on $\lambda$, hence a critical behavior may just be expected at $\rho_1$ in the limit
$\lambda \to \rho_1 g'_{\infty}$.\\

Summarizing these observation we obtain: if
\[
\lim_{\lambda \to \rho_1 g'_{\infty}}  \int_{\rho_1}^{\rho_2} (g')^{-1}(\lambda/r) \dr = \infty \psp ,
\]
then $\Delta m(\lambda)$ takes any value in $(0,\infty)$ and for all $m_1 < m_2 \in \rz$
problem \reff{intro 11} admits a solution taking the boundary data.\\

If
\[
\lim_{\lambda \to \rho_1 g'_{\infty}} \int_{\rho_1}^{\rho_2} (g')^{-1}(\lambda/r) \dr
=: \Delta m_\infty < \infty \psp ,
\]
then a solution taking the boundary data exists if and only if $m_2 -m_1 < \Delta m_\infty$.\\

At this point we note that, given $\zeta_1 < \zeta_2$ such that $\uh_{\zeta_1}(\rho_1)= \uh_{\zeta_2}(\rho_2)$,
the monotonicity immediately shows $\uh_{\zeta_1}(\rho_1) < \uh_{\zeta_2}(\rho_1)$
(compare \reff{mono 1} and \reff{mono 2} in the case $\zeta_2\omax \leq \zeta_1\omax$). \\

Let us finally suppose that $\zeta_1 < \zeta_2 < m_2 - \Delta m_\infty$ for
some real number $\zeta$. By the above considerations we have
\[
\uh_{\zeta_1}(\rho_1) = \uh_{\zeta_2}(\rho_1) = m_2 - \Delta m_\infty \psp ,
\]
and the limit number $m_2 - \Delta m_\infty$ serves as the boundary datum
for $\uh_{\zeta_1}$ as well as for $\uh_{\zeta_2}$, which immediately gives the corollary.
\hfill\qed\\

\section{Examples}\label{ex}

We finally sketch three characteristic examples by presenting explicit solutions.\\

To this purpose we recall the one parameter family given in \reff{intro 6} (now denoted by $g_\mu$)
\[
g_\mu(t) := \left\{\begin{array}{rcl}
\dis t - \frac{1}{2-\mu} (1 + t)^{2- \mu} -  \frac{1}{\mu - 2}&\mbox{if}& \mu\not=2\psp ,\\[2ex]
\dis t - \ln(1+t)&\mbox{if}& \mu =2 \psp .
\end{array}\right.
\]

Note that ${g_\mu}'_\infty = 1$ for any $\mu >1$.\\

With this choice of $g_\mu$, the Euler equation \reff{cor 5} reads as
\begin{equation}\label{ex 2}
\hat{u}' = \Bigg[\frac{r}{r-\lambda}\Bigg]^{\frac{1}{\mu-1}} - 1 \psp .
\end{equation}

We note that the condition (1.9) of \cite{BBM:2018_1} motivates to consider 
examples choosing $1<\mu < 2$, $\mu=2$, $\mu > 3$, respectively. 
\begin{enumerate}
\item Suppose that $\mu =3/2$, i.e.~$1< \mu < 2$. Then 
\[
\hat{u}(r) = 2 \lambda \ln(r-\lambda) - \frac{\lambda^{2}}{r-\lambda}+c  \psp , \fsp c\in \rz\psp ,
\]
provides an exact solution. With the notation from above we have
\[
\Delta m(\lambda)  := \Bigg[2 \lambda \ln (\rho_2 -\lambda) - \frac{\lambda^2}{\rho_2-\lambda}\Bigg]
-\Bigg[2 \lambda \ln (\rho_1 -\lambda) - \frac{\lambda^2}{\rho_1-\lambda}\Bigg] \psp .
\]

We see that in this case we have
\[
\Delta m(\lambda)\to \infty\fsp\mbox{as}\msp \lambda \to \rho_1 \psp ,
\] 
and with the right choice of the free parameters $\lambda$, $c$  we obtain a smooth solution to 
\reff{intro 11} taking the boundary data.

\item Consider the limit case $\mu =2$. We then have
\[
\hat{u}(r) = \lambda \ln(r-\lambda) + c
\]
and as above we find for any given boundary data a solution realizing this data.

\item In the case $\mu =3$ we find as solution of the Euler equation
\[
\hat{u}(r) = \sqrt{r^2 - r \lambda}- t + 
\frac{\lambda}{2} \ln\Bigg[\frac{2r-\lambda + 2 \sqrt{r^2 - r \lambda}}{2}+ c\Bigg]\psp .
\]
Now we note that 
\[
\Delta m_\infty < \infty \psp ,
\]
hence that boundary data can not be attained if $\Delta m_\infty < m_2 - m_1$.
\end{enumerate}

\bibliography{radial}
\bibliographystyle{unsrt}

\begin{tabular}{ll}
Michael Bildhauer&bibi@math.uni-sb.de\\
Martin Fuchs&fuchs@math.uni-sb.de\\[5ex]
Department of Mathematics&\\
Saarland University&\\
P.O.~Box 15 11 50&\\
66041 Saarbr\"ucken&\\ 
Germany&
\end{tabular}

\end{document}